\documentclass{conm-p-l}

\newtheorem{theorem}{Theorem}[section]

\newtheorem{proposition}{Proposition}[section]
\newtheorem{corollary}{Corollary}[section]
\theoremstyle{definition}
\newtheorem{definition}[theorem]{Definition}

\theoremstyle{remark}
\newtheorem{remark}[theorem]{Remark}

\numberwithin{equation}{section}

\newcommand{\CC}{\mathbb{C}}   

\newcommand{\HH}{\mathbb{H}}   

\newcommand{\RR}{\mathbb{R}}   

\newcommand{\ZZ}{\mathbb{Z}}

\newcommand{\Ca}{\mathbb{C}\textrm{a}}

\newcommand{\II}{\textrm{Im}}

\begin{document}

\title[Cayley 4-Frames and a Quaternion K\"ahler Reduction]{Cayley 4-Frames 
and\\ a Quaternion K\"ahler Reduction Related to 
$\boldsymbol{\mathrm{Spin}(7)}$}

\author{Liviu Ornea}
\address{University of Bucharest, Faculty of Mathematics, 14 Academiei str., 70109
Bu\-cha\-rest, Romania}
\email{lornea@imar.ro}

\author{Paolo Piccinni}
\address{Dipartimento di Matematica, Universit\`a degli Studi di Roma "La Sapienza", 
Piazzale Aldo Moro 2, I-00185, Roma, Italy}
\email{piccinni@mat.uniroma1.it}

\thanks{The authors are members of EDGE, Research Training Network HRPN-CT-2000-00101, supported by The European Human Potential Programme.}



\begin{abstract} The object of this note is CAYLEY, the Grassmannian of the oriented 4-planes in $\RR^8$ that are
closed under the three-fold cross product. We describe an action of $\mathrm{U}(1) \times
\mathrm{Sp}(1)$ on the quaternionic projective space $\HH P^7$, that allows to obtain a
$\ZZ_2$-quotient of CAYLEY by quaternion K\"ahler reduction.
\end{abstract}

\maketitle

\section{Introduction}

The existence of only two exceptional cross products - in $\RR^7$ and in $\RR^8$, with two and three factors respectively -
attracted the interest of Alfred Gray in the sixties, and this was one of his approaches to the study of holonomies
$\mathrm{G_2}$ and $\mathrm{Spin}(7)$ on Riemannian manifolds \cite{br-gr}, \cite{Gr2}. About one decade later, the symmetric
space structure of the Grassmannians of those planes in $\RR^7$ or $\RR^8$ that are closed under such cross products was
recognized in the positive quaternion K\"ahler manifolds
$\frac{\mathrm{G_2}}{\mathrm{SO}(4)}$ and $\frac{\mathrm{Spin}(7)}{(\mathrm{Sp}(1)
\times \mathrm{Sp}(1) \times \mathrm{Sp}(1))/\ZZ_2}$ \cite{mar}, \cite{ha-la}. The latter of these manifolds is in fact
isometric to the Grassmannian ${\mathrm{Gr}}_4(\RR^7) =
\frac{\mathrm{SO}(7)}{\mathrm{SO}(3) \times \mathrm{SO}(4)}$ of oriented 4-planes in $\RR^7$, but its r\^ole as "exceptional
Grassmannian" of some distinguished 4-planes in $\RR^8$ is so interesting to make it deserving of notations like CAY or
CAYLEY in papers on calibrations \cite{br-ha},
\cite{gl-mc-mo}.

In this note we describe how CAYLEY can be obtained -- up to a $\ZZ_2$-quotient -- through a quaternion K\"ahler reduction of
the projective space $\HH P^7$, acted on by a group isomorphic to $\mathrm{U}(1) \times \mathrm{Sp}(1)$. This action is
similar to that used by P. Kobak and A. Swann in $\HH P^6$, obtaining a $\ZZ_3$-quotient of
$\frac{\mathrm{G_2}}{\mathrm{SO}(4)}$ by quaternion K\"ahler reduction
\cite{ko-sw}. In the later note \cite{ko-sw2} a different action of the same group in $\HH P^7$ is
described, obtaining this time as a reduction a $\ZZ_2$-quotient of ${\mathrm{Gr}}_4(\RR^7)$. Our reduction is in fact
equivalent to this latter, via the isometry CAYLEY $\cong {\mathrm{Gr}}_4(\RR^7)$, although our definition of the
action is much closer to the point of view of the former article
\cite{ko-sw}. We observe also that, through the same isometry with this real Grassmannian, one can obtain CAYLEY as a
reduction of $\HH P^6$ acted on by
$\mathrm{Sp}(1)$, following one of the classical procedures described in \cite{Bo-Ga} (and in this way no finite quotient
is involved). However, the reduction we are going to describe here -- without making use of the isometry  CAYLEY $\cong
{\mathrm{Gr}}_4(\RR^7)$ -- has the advantage of admitting two interesting generalizations. 

One of them is a consequence of the possibility of introducing weights in the action of the
$\mathrm{U}(1)$-factor. This allows to obtain a family of 12-dimensional quaternion K\"ahler orbifolds, some
of them admitting smooth 3-Sasakian manifolds over them. These smooth 15-dimensional manifolds are, together with
similar 11-dimensional manifolds related to $\mathrm{G_2}$, the first examples of 3-Sasakian manifolds which are neither
homogeneous nor toric \cite{bo-ga-pi}. 

The other possibility of extending the present reduction is obtained by looking at
two other quaternion K\"ahler Wolf spaces. Since CAYLEY can be regarded
as the manifold of the hypercomplex 4-planes in $ \RR^8$ identified with the real vector space of Cayley numbers $\Ca$ (Proposition 3.1 below), higher
dimensional analogues of it are the manifolds $\frac{\mathrm{Spin}(9)}{(\mathrm{Sp}(2) \times \mathrm{Sp}(1) \times
\mathrm{Sp}(1))/\ZZ_2}
\cong {\mathrm{Gr}}_4(\RR^9)$ and the exceptional Wolf space $\frac{\mathrm{F_4}}{\mathrm{Sp}(3)
\cdot \mathrm{Sp}(1)}$, geometrically the manifolds of the $\HH P^1 \subset \Ca P^1$ and of the $\HH P^2 \subset \Ca P^2$,
respectively. Some issues related to this second generalization will be studied in a future work \cite{BGOP}.

{\bf Acknowledgements.~}\;The first author acknowledges financial support by C.N.R. 
of Italy and by the Cultural Agreement between Universit\`a di Roma "La Sapienza" and
University of Bucharest. Both authors thank Kris Galicki for helpful and stimulating conversations and Robert Bryant
for a kind observation.

\section{Preliminaries on Cayley numbers}

Let $\Ca$ be the algebra of Cayley numbers, and let 
$\{1,{\mathrm{i,j,k,e,f=ie,g=je,}}$ ${\mathrm{h=ke}}\}$ be its canonical basis over $\RR$.
The multiplication is given by $$xy=(ac-\overline{d} b )+(b \overline{c}+ d a)\mathrm e,$$ where $x=a+b\mathrm e, y=c+d\mathrm e \in
\Ca$  are written through the identification $\Ca \cong \HH^2$ with pairs of quaternions.
The quaternionic conjugation (already used in the previous formula) induces a conjugation in
$\Ca$: $\overline{x}=\overline a - b\mathrm e$, allowing to write the non-commutativity rule: $\overline{xy}=
{\overline y} ~{\overline x}.$ 

The non-associativity of $\Ca$ gives rise to the associator $[x,y,z]=(xy)z-x(yz),$
alternating form that vanishes whenever two of its arguments are either equal
or conjugate. Geometrically, the associator defines the class of \emph{associative 3-planes} in
$\RR^7 \cong \II \Ca$, defined in orthonormal bases by $[x,y,z]=0$. They are characterized as the
3-planes of $\RR^7$ closed with respect to the \emph{two-fold cross product}: 
$$x\times y = \II (\overline y x),$$ that is so more generally defined for any $x,y \in \RR^8 \cong \Ca$. Note that
if $x,y$ are orthogonal and in $\RR^7 \cong \II \Ca$, the cross product is simply $xy$.  

The \emph{three-fold cross product} in $\Ca\cong \RR^8$ is defined by the formula $$x\times y\times
z=\frac{1}{2}(x({\overline y}z)-z({\overline y}x)),$$
reducing to $x(\overline y z)$ for $x,y,z$ orthogonal.

The following properties hold whenever $x,y$ are orthogonal and for any $w \in \Ca$: 
\begin{equation}
x(\overline y w) = - y(\overline x w), \qquad (w \overline y)x = -(w\overline x)y.
\end{equation}
Moreover, for any $x,y,z \in \Ca$:
\begin{equation}
(xy)(zx) = x(yz)x.
\end{equation}
(A reference for these preliminaries is \cite{ha-la}, Appendix IV).

\section{The Stiefel manifold of Cayley 4-frames}

A 4-plane $\zeta$ of $\RR^8$ that is closed under the three-fold 
cross product is called a \emph{Cayley 4-plane} and it is
oriented by choices of bases $\{w=x\times y \times z,x,y,z\}$. The manifold of the
Cayley 4-planes in $\RR^8$ is $\text{CAYLEY}~Ê=~\frac{\mathrm{Spin}(7)}{(\mathrm{Sp}(1) \times \mathrm{Sp}(1) \times
\mathrm{Sp}(1)/\ZZ_2)},$  12-dimensional quaternionic submanifold of the Grassmannian ${\mathrm{Gr}}_4(\RR^8)$ of oriented 4-planes in
$\RR^8$ (\cite{mar}, p. 262 or \cite{ha-la}, p. 123).

Proposition IV,1.27 in \cite{ha-la} states that a 4-plane $\zeta$ in $\RR^8$ is Cayley if and only if $-\zeta$ is closed
under the complex structures defined by the 2-planes $\alpha \subset \zeta$.  This fact can be reformulated as follows.  
\begin{proposition}
A 4-plane $\zeta$ in $\RR^8$ is Cayley if and only if any triple of mutually orthogonal 2-planes $\alpha, \beta, \gamma
\subset \zeta$, all intersecting in a line, defines a hypercomplex structure on $\zeta$.
\end{proposition}
\begin{proof} 
For any $\zeta$,
$\dim (\zeta\cap \II \Ca)$ is either  3 or 4. Thus, if $\zeta \in$ CAYLEY we may select
orthonormal imaginary octonions $x,y,z \in \zeta$ such that
$\{x\times y\times z,x,y,z\}$ is an oriented basis of $\zeta$. If $u=x\times (x\times y\times z)=y\times z$,
$v=y \times (x\times y\times z)=z\times x$, $w=z\times (x\times y\times z)=x \times y$ we have $u, v, w\in S^6$, and their
corresponding complex structures $J_u,J_v,J_w$ are associated to the 2-planes $\alpha =$ span$\{x\times y \times z,
x\}$, $\beta =$ span$\{x\times y \times z, y\}$,$\gamma =$ span$\{x\times y \times z, z\}$. Since $J_u y = -z$, $J_v x =
z$, $J_w x = -y$, then $(J_u,J_v,J_w)$ satisfy 
$J_v \circ J_u =- J_u \circ J_v = J_w$, i.e. it is a hypercomplex structure
on
$\zeta$. The converse follows from the aformentioned characterization in \cite{ha-la}, p. 119. 
\end{proof}

Our construction of $(u,v,w)$ out of $\zeta=$ span $\{x \times y \times z, x, y, z\}$
corresponds to the isometry $\sim:$ CAYLEY $\rightarrow {\mathrm{Gr}}_3(\RR^7)$ of \cite{br-ha}, p. 11. The image
under $\sim$ of the $\zeta \in$ CAYLEY can be interpreted as a \emph{tricomplex section} of $S^6$, oriented orthonormal
bases of the $\zeta^\sim$ being triples $u,v,w$ of unit octonions non necessarily satisfying the hypercomplex
relations. The non-associativity of $\Ca$ allows and ensures that such triples define a hypercomplex structure on $\zeta$.
An example is the Cayley 4-plane
$\zeta$ = span $\{1-{\mathrm{h,i+g,j-f,k+e}}\}$: our procedure 
gives the tricomplex triple $(u,v,w)=({\mathrm{i,j,e}})$, whose associated
$(J_{\mathrm i},J_{\mathrm j},J_{\mathrm e})$ is hypercomplex on
$\zeta$. 

This discussion permits to describe the inverse of the isometry $\sim$, as follows.
\begin{corollary} 
Given a tricomplex section of $S^6$ with oriented orthonormal basis 
$(u,v,w)$, there is a unique Cayley 4-plane $\zeta$ in $\RR^8$
on which $(J_u,J_v,J_w)$ is hypercomplex.
\end{corollary}

\begin{definition}
A \emph{Cayley 4-frame} in $\RR^8$ is an oriented orthonormal 4-frame in a Cayley 4-plane $\zeta$, hence a frame $\{x,
I_1x, I_2x, I_3x\}$, where
$(I_1,I_2,I_3)$ is the hypercomplex structure of $\zeta$.
\end{definition}
By the action of
$\mathrm{Spin}(7) \supset \mathrm{G_2} \supset \mathrm{SU}(3)$ on the spheres $S^7 \supset S^6 \supset S^5$, the
latter with isotropy $\mathrm{SU}(2) \cong \mathrm{Sp}(1)$, we have:
\begin{proposition}\label{cay}
The Stiefel manifold of Cayley 4-frames in $\RR^8$ is the homogeneous space
$V=\frac{\mathrm{Spin}(7)}{\mathrm{Sp}(1)}$.
\end{proposition}

Observe finally that:
\begin{proposition}\label{doi}
An orthonormal frame $\{f_1,f_2,f_3,f_4\}$ in $\RR^8$ is a Cayley 4-frame if and only if $\overline{f}_2 f_1
=\overline{f}_3 f_4$.
\end{proposition}

\section{CAYLEY and a reduction of $\HH P^7$.}

We now show how a $\ZZ_2$-quotient of CAYLEY can be obtained as 
quaternion K\"ahler reduction of  $\HH P^7$ by the action of $\mathrm{U}(1)\times \mathrm{Sp}(1)$.  
According to \cite{Bo-Ga-Ma 1}, we first  reduce (by the same group) the 3-Sasakian manifold 
which stands over $\HH P^7$, namely the sphere $S^{31}$. Then we interpret  the quotient of 
$S^{31}$ by $\mathrm{U}(1)\times \mathrm{Sp}(1)$ as the total space of an $\mathrm{SO}(3)$-bundle over a 
quaternion K\"ahler orbifold which is the quotient of 
 $\HH P^7$ by the same group. 

The 3-Sasakian sphere $S^{31}$ is acted on by $\mathrm{U}(1)\times
\mathrm{Sp}(1)$ as follows. The factor $\mathrm{Sp}(1)$ acts by right multiplication on ${\vec h} =
(h_\alpha) \in S^{31}\subset
\HH^8$, and the moment map  $\mu:S^{31}\rightarrow \RR^9$ of the action reads: 
$$\mu (\vec h)=(\sum_{\alpha=1}^8\overline h_\alpha \mathrm i h_\alpha,\quad
\sum_{\alpha=1}^8\overline h_\alpha \mathrm j h_\alpha,\quad
\sum_{\alpha=1}^8\overline h_\alpha \mathrm k h_\alpha).$$ By writing $\vec h =\vec a +\vec b \mathrm i+\vec c \mathrm j+ \vec d
\mathrm k$,  it is easy to see that $\mu^{-1}(0)$ coincides with the Stiefel manifold of oriented (renormalized) orthonormal 4-frames
in
$\RR^8$ \cite{Bo-Ga-Ma 1}. 

We then act by the factor $\mathrm{U}(1)$, rotating pairs of 
coordinates. This is explicitly described by $\vec h \mapsto
\mathrm{diag}\left(A(\theta),A(\theta), A(\theta),A(\theta)\right)\cdot \vec h$, where
$A(\theta)=\left(\begin{smallmatrix}\cos
\theta& -\sin
\theta\\
\sin
\theta&\cos
\theta\end{smallmatrix}\right)$, $\theta\in \RR$.
The
associated moment map
$\nu: S^{31}\rightarrow 
\RR^3$ is now:
$$\nu (\vec h)= \sum_{\beta =1}^4(\overline h_{2\beta -1}h_{2\beta}-\overline
h_{2\beta}h_{2\beta -1}).$$ We are interested in the common zero set $\mathcal{N}=\mu^{-1}(0) \cap \nu^{-1}(0)$. 
\begin{proposition}
$\mathcal{N}=\mathrm{U}(1) \cdot V $, where $\cdot$ is the action of
$\mathrm{U}(1)$ on Cayley 4-frames.
\end{proposition}
\begin{proof}
The inclusion $V \subset \mathcal N$ can be checked either by direct computation, using Proposition 3.3, or by a standard
choice of the frame, like $(1,{\mathrm{i,j,k}})$, and the observation 
that  $\nu(\vec h) =\nu(\vec a +\vec b {\mathrm i} + \vec c {\mathrm j} 
+ \vec d {\mathrm k}) =0$
is invariant under right multiplication of $\vec a, \vec b, \vec c, \vec d$ by any $u \in S^6$, and hence by $\mathrm{Spin}(7)$
(cf. \cite{ha-la}, p. 121). It follows also
$\mathrm{U}(1) \cdot V \subset \mathcal N$, by the $\mathrm{U}(1)$-equivariance of $\nu$.
\par Conversely, to see that $\mathcal N \subset \mathrm{U}(1) \cdot V $, refer to a standard choice of three vectors 
to be substituted in the moment map equation $\nu(\vec h)  = 
\nu(\vec a +\vec b {\mathrm i}
+ \vec c {\mathrm j} + \vec d {\mathrm k})=0$, assuming $f_2 = \vec b ={\mathrm j},
f_3 = \vec c ={\mathrm e}, f_4 =\vec d= {\mathrm g}$,  (cf. the similar proof of the
$\mathrm{G_2}$-case in
\cite{ko-sw}). Then the equation $\nu(\vec h) =0$ and the orthonormality of the frame give $\vec f_1 = \vec a = \cos \theta +
\sin \theta {\mathrm i}$. Then it is easy to check that the element ${\mathrm e}^{-\mathrm i\frac{\theta}{2}}$ of $\mathrm{U}(1)$
transforms
$(\cos \theta +
\sin \theta {\mathrm i},{\mathrm{j,e,g}})$ into a Cayley 4-frame.
\end{proof}

Observe now that $\mathrm{U}(1) \cap \mathrm{Spin}(7) = \mathrm{U}(1) \cap \mathrm{SU}(4) = \ZZ_4$ with generator
$\tau=\mathrm e^{\mathrm i\frac{\pi}{2}}$, and that under the action of $\tau$ on $V$, a Cayley 4-frame $(f_1,f_2,f_3,f_4)$ is
transformed into another frame of the same Cayley 4-plane if and only if $(f_1,f_2,f_3,f_4)$ is complex unitary, i.e. an element of 
$\frac{\mathrm{SU}(4)}{\mathrm{Sp}(1)}$. Also, of
course
$\tau^{2} = -1$, so that any Cayley 4-plane is fixed under it. This explains the
following description of the orbits of the $\mathrm{U}(1) \times \mathrm{Sp}(1)$-action on
$\mathcal{N}$: points in $\frac{\mathrm{SU}(4)}{\mathrm{Sp}(1)}\subset V$ generate orbits that are the fixed points of an
induced action of $\ZZ_2$ on all the orbits of $\mathcal{N}$, and a 3-Sasakian orbifold $\ZZ_2 \backslash
\frac{\mathrm{Spin}(7)}{\mathrm{Spin}(4)}$ is obtained as quotient. We state the corresponding quaternion
K\"ahler reduction. 
\vspace{-.1in}
\begin{theorem}
The quaternion K\"ahler quotient of $\HH P^7$ by the described action of $\mathrm{U}(1) \times \mathrm{Sp}(1)$ is 
an orbifold $\ZZ_2\backslash \text{\rm CAYLEY}$, with a singular stratum isometric to the complex Grassmannian
$\frac{\mathrm{SU}(4)}{S(\mathrm{U}(2) \times
\mathrm{U}(2))}$.
\end{theorem}
\vspace{-.1in}
\begin{remark} By identifying any $\zeta$ with its
orthogonal complement $\zeta^\perp$, one obtains a \emph{smooth} $\ZZ_2$-quotient of CAYLEY. Since $\perp$
corresponds to the change of orientation on 4-planes in
$\RR^7$ (\cite{br-ha}, p. 11), this smooth $\ZZ_2$-quotient of CAYLEY is the locally
quaternion K\"ahler Grassmannian of \emph{unoriented} 4-planes in $\RR^7$. 
\par The $\ZZ_2\backslash$ CAYLEY given by Theorem 4.1 is not smooth, its
construction yielding the stratified space $\mathcal M_{\text{reg}} \cup
$ ${\mathrm{Gr}}_2(\CC^4)$. The singular stratum ${\mathrm{Gr}}_2(\CC^4)$ corresponds, under the isometry CAYLEY
$\cong {\mathrm{Gr}}_4(\RR^7)$, to the standard
${\mathrm{Gr}}_4(\RR^6) \subset {\mathrm{Gr}}_4(\RR^7)$. Thus the orbifold
$\ZZ_2\backslash \text{\rm CAYLEY}$ in Theorem 4.1 is isometric to the singular quotient
${\mathrm{Gr}}_4(\RR^7)/\sigma_{\RR^6}$ by the symmetry $\sigma_{\RR^6}$ with respect to $\RR^6 \subset \RR^7$.    
\end{remark}

\end{document}